\documentclass{article}

\usepackage{microtype}
\usepackage{graphicx}
\usepackage{subfigure}
\usepackage{booktabs} 

\usepackage{hyperref}

\usepackage{multirow}

\usepackage[accepted]{icml2023}

\usepackage{amsmath}
\usepackage{amssymb}
\usepackage{mathtools}
\usepackage{amsthm}
\usepackage{url}

\usepackage[capitalize,noabbrev]{cleveref}

\theoremstyle{plain}

\theoremstyle{definition}

\theoremstyle{remark}

\DeclareMathOperator{\grad}{\nabla}
\DeclareMathOperator*{\argmin}{arg\,min}
\newcommand{\RR}{\mathbb{R}}
\newcommand{\eps}{\varepsilon}
\newcommand{\norm}[1]{\left\lVert #1 \right\rVert}
\newcommand{\abs}[1]{\left\lvert #1 \right\rvert}
\newcommand{\D}{\mathrm{d}}

\usepackage{listings}       

\definecolor{codeblue}{rgb}{0,0,1}
\definecolor{codegrey}{rgb}{0.5,0.5,0.5}
\definecolor{codepurple}{rgb}{0.58,0,0.82}
\definecolor{backcolor}{rgb}{0.95, 0.95, 0.95}
\definecolor{stringcolor}{rgb}{0,0.63,0} 

\lstdefinelanguage{Python}{ 
numbers=left, 
numberstyle=\tiny, 
numbersep=1em, 
xleftmargin=1em, 
framextopmargin=2em, 
framexbottommargin=2em, 
showspaces=false, 
showtabs=false, 
showstringspaces=false, 
frame=l, 
tabsize=4, 
basicstyle=\small\ttfamily,
backgroundcolor=\color{backcolor}, 
commentstyle=\color{codegrey}\slshape, 
stringstyle=\color{stringcolor}, 
morecomment=[s][\color{stringcolor}]{"""}{"""}, 
morecomment=[s][\color{stringcolor}]{'''}{'''}, 
morekeywords={import,from,class,def,for,while,if,is,in,elif,else,not,and,or,print,break,continue,return,True,False,None,access,as,,del,except,exec,finally,global,import,lambda,print,raise,try,assert}, 
keywordstyle={\color{codeblue}\small\ttfamily}, 
keywordstyle={[2]\color{codeblue}\slshape}, 
emphstyle={\color{self}\slshape}, 
}  
\lstnewenvironment{python}[1][]{%
  \lstset{language=Python,#1}%
}{}

\usepackage[textsize=tiny]{todonotes}

\icmltitlerunning{Optimistix}

\begin{document}

\twocolumn[
\icmltitle{Optimistix: modular optimisation in JAX and Equinox}




\begin{icmlauthorlist}
\icmlauthor{Jason Rader}{oxford}
\icmlauthor{Terry Lyons}{oxford}
\icmlauthor{Patrick Kidger}{googlex}
\end{icmlauthorlist}

\icmlaffiliation{oxford}{Mathematical Institute, University of Oxford}
\icmlaffiliation{googlex}{Google X}

\icmlcorrespondingauthor{Jason Rader}{rader@maths.ox.ac.uk}
\icmlcorrespondingauthor{Patrick Kidger}{math@kidger.site}

\icmlkeywords{Machine Learning}

\vskip 0.3in
]



\printAffiliationsAndNotice{} 

\begin{abstract}
We introduce Optimistix: a nonlinear optimisation library built in JAX and Equinox. Optimistix introduces a novel, modular approach for its minimisers and least-squares solvers. This modularity relies on new practical abstractions for optimisation which we call search and descent, and which generalise classical notions of line search, trust-region, and learning-rate algorithms. It provides high-level APIs and solvers for minimisation, nonlinear least-squares, root-finding, and fixed-point iteration. Optimistix is available at  https://github.com/patrick-kidger/optimistix.
\end{abstract}
\section{Introduction}
JAX is a Python autodifferentiation framework popular for scientific computing and machine learning \cite{Bradbury2018, KidgerThesis}. Equinox \cite{Kidger2021} extends many JAX core transformations and concepts, and adds additional functionality for parameterised functions. Equinox has become a popular choice for machine learning \cite{Haliax, Levanter, Eqxvision, palmjax} and scientific machine learning `sciML' \cite{KidgerThesis, lineax, pastrana_jaxfdm_2023, flowmc} in JAX.

We introduce Optimistix, a nonlinear optimisation library built in JAX + Equinox. Optimistix targets differentiable scientific computing and sciML tasks. The sciML ecosystem in JAX is large and growing, and already includes packages for differentiable rigid-body physics simulation \cite{brax2021github}, computational fluid dynamics \cite{Dresdner2022-Spectral-ML, BEZGIN2022108527}, protein structure prediction \cite{AlphaFold2021}, linear solves and least-squares \cite{lineax}, ordinary and stochastic differential equations \cite{KidgerThesis}, general-purpose optimisation \cite{jaxopt_implicit_diff}, structural design \cite{pastrana_jaxfdm_2023}, Bayesian optimisation \cite{oss_vizier, google_vizier}, and probabilistic modeling \cite{synjax2023}.

\subsection{Contributions}
\label{sec:contributions}
We introduce Optimistix, a JAX nonlinear optimisation library with the following features:
\begin{enumerate}
    \itemsep0em 
    \item Modular optimisers.
    \item Fast compile times and run times.
    \item Support for general PyTrees, \footnote{JAX data structures consisting of arbitrarily nested container types, containing other JAX/Python types. The containers (tuples/dictionaries/lists/custom types) are referred to as `nodes', and the data types they hold as `leaves'.} and use of PyTrees for solver state.
\end{enumerate}
We highlight the first point as the central contribution of this work. This allows users to define custom optimisers for a specific problem by swapping components of the optimiser.

To achieve this, Optimistix introduces two practical abstractions: search and descent. These abstractions generalise classical line search, trust-region, and learning-rate algorithms. To the best of our knowledge, Optimistix is the first optimisation software modularised in this way.

Optimistix includes APIs for four types of optimisation tasks: minimisation, nonlinear least-squares, root-finding, and fixed-point iteration. Each has a high-level API, and support automatic conversion of appropriate problem types.

Optimistix is already seeing adoption in sciML. For example, in software libraries for solving differential equations \cite{KidgerThesis} and probabilistic inference \cite{bayeux}.

Example usage on the Rosenbrock problem \cite{rosenbrock} using BFGS and automatically converting least-squares to minimisation is:
\begin{lstlisting}[language=python]
import jax
import optimistix as optx

def loss(x, scaling):
  res1 = scaling * (x[1:] - x[:-1]**2)
  return (res1, 1 - x)

init = jax.numpy.zeros(100)
scaling = 10
solver = optx.BFGS(rtol=1e-5, atol=1e-6)
minimum = optx.least_squares(
  loss, solver, init, args=scaling
)
\end{lstlisting}

    


\section{Background: Preconditioned Gradient Methods}
\label{sec:quadratic}

The dominant approach to differentiable optimisation is to locally approximate the objective function $f: \RR^N \to \RR$ at an iterate $x_k$ using a quadratic model function:
\begin{equation}
    m_k(p) = f(x_k) +  \nabla f(x_k)^T p + \frac{1}{2} p^T H_k p \label{eqn:preconditioned}
\end{equation}
\noindent
where $p \in \RR^N$, and $H_k \in \RR^{N \times N} $ is an approximation to the Hessian $\nabla^2 f(x_k)$ \cite{Nocedal2006}[sections 1-6], \cite{Bonnans2006}[section 4], \cite{Conn2000}. This quadratic is used to find the next step in an iterative algorithm to minimise $f$. The minimum of (\ref{eqn:preconditioned}) is found at $p = -H_k^{-1} \nabla f(x_k)$. As such, these are sometimes called preconditioned gradient methods \cite{gupta2018shampoo} (not to be confused with preconditioned conjugate gradient methods \cite{Nocedal2006}[section 5].)

Most first-order methods common in the machine learning literature, such as Adam \cite{Kingma2017} and Adagrad \cite{Duchi2011}, are also preconditioned gradient methods. In these algorithms, $H^{-1}$ is usually stored and updated directly, rather than computed from $H$. Moreover, $H^{-1}$ is often highly structured or diagonal to reduce memory cost.


The class of preconditioned gradient methods is extremely large, and includes Newton's method, quasi-Newton methods, and Gauss-Newton methods \cite{Bonnans2006}[section 4], \cite{Nocedal2006}[section 6], some nonlinear conjugate gradient methods \cite{Sherali1990}, gradient descent, and adaptive gradient methods such as Adam \cite{Kingma2017} and Adagrad \cite{Duchi2011}. In many of these algorithms, the quadratic approximation is implicit.


Preconditioned gradient methods are central in Optimistix, and were an important motivator of the descent and search abstractions we now introduce in section \ref{sec:modularity}.



\section{Modularity in Optimistix} \label{sec:modularity}

This section introduces the key advancement of Optimistix: modularity. To attain this, Optimistix uses a generalised approach to line searches, trust-regions, and learning-rates through the abstractions of search and descent. These concepts offer a precise formalism of ideas already present in the optimisation literature. This approach makes it easier both to advance theory (in Section \ref{sec:mixing} we demonstrate how this may be used to create a novel minimiser), and serve as a practical approach to modularity. To the best of our knowledge this approach is new here, and is not present in any other open-source optimisation package.

Consider a scalar function to minimise: $f:\RR^n \to \RR$. Searches consume local information about $f$, such as its value, gradient, and/or Hessian at a point, and return a scalar. This scalar corresponds to the distance along a line search, the trust-region radius of a trust-region method, the value of a learning-rate, etc. Searches are a generalisation of these.

A descent consumes this scalar, as well as the same local information about $f$, and returns the step the optimiser should take. For example, gradient descent with a fixed learning-rate $\alpha$ is implemented as a search which returns a fixed value $\alpha \in \RR$, and a descent which returns $\alpha \grad f(x_k)$.

Searches and descents are modular components in Optimistix, and they can be easily swapped by the user. 

\subsection{Creating a Novel Optimiser with Ease: an Example} \label{sec:mixing}

Optimistix provides a new "mix-and-match" API which allows users to easily create new optimisers. For example, here we define a custom optimiser and use this to solve a toy linear regression problem:
\begin{lstlisting}[language=Python]
from collections.abc import Callable
import jax.numpy as jnp
import jax.random as jr
import optimistix as optx

class HybridMinimiser(optx.AbstractBFGS):
  rtol: float
  atol: float
  norm: Callable = optx.max_norm
  use_inverse: bool = False
  descent: optx.AbstractDescent = (
        optx.DoglegDescent()
    )
  search: optx.AbstractSearch = (
        optx.LearningRate(0.1)
    )

solver = HybridMinimiser(
  rtol=1e-8, atol=1e-9
)

def loss(weights, data):
  x, y = data
  return weights.T @ x - y

key1, key2 = jr.split(jr.PRNGKey(0))
noise = 0.1 * jr.normal(key1, shape=(99,))

true_weights = jnp.array([3.14, -7, 2.71])
xs = jr.normal(key2, shape=(3, 99))
ys = true_weights.T @ xs + noise
init = jnp.array([0.0, 0.0, 0.0])

# Convert least-squares to minimisation
# then solve.
minimum = optx.least_squares(
  loss, solver, init, args=(xs, ys)
)
# `minimum.value`:
# Array([3.129, -6.983,  2.732]
\end{lstlisting}
\newpage
\noindent
\verb|HybridMinimiser| defines an optimiser which, at each step $k$, uses the gradient and BFGS quasi-Newton approximation $B(x_k) \approx \nabla^2 f(x_k)$ \cite{Nocedal2006}[section 6.1] to construct a quadratic model function as discussed in section \ref{sec:quadratic}. A `dogleg' descent path is then built, which interpolates between $x_k$, $\grad f(x_k)$, and $B^{-1}(x_k) \grad f(x_k)$ with a piecewise linear curve. Finally, a constant step of length \verb|0.1| is taken along this descent path, to form the next iterate $x_{k + 1}$.

\verb|HybridMinimiser| is not an off-the-shelf optimiser, and would require a custom implementation in most optimisation packages. Implementing a custom, performant algorithm may require significant experience in optimisation, and hundreds of lines of technical code. In Optimistix, a performant implementation of this novel optimiser takes fewer than 10 lines of code.

Custom optimisers allow users to choose optimisation methods appropriate for the problem at hand. For example, the novel optimiser above can solve the poorly-scaled Biggs EXP6 function \cite{tuos}, whereas the standard BFGS algorithm with backtracking line search fails to solve the problem to an acceptable accuracy.



\subsection{Search}

For most nonlinear functions, the quadratic approximation $m_k$ in (\ref{eqn:preconditioned}) is only a reasonable approximation in a small neighborhood of $x_k$. The full preconditioned gradient step $-H_k^{-1} \grad f(x_k)$ often overshoots the region where the approximation is good, slowing the optimisation process. Line searches, trust-regions, and learning-rates are all methods to keep steps within a region where the approximation is good.

We begin by discussing the classical approach to line searches and trust-regions.

\textbf{Line searches}

Line searches \cite{Nocedal2006}[section 3] move in the direction of the preconditioned gradient, but only move a certain amount $\alpha_k \in \RR^+$. ie.
\begin{equation*}
    x_{k + 1} = x_k - \alpha_k H_k^{-1} \grad f(x_k).
\end{equation*}
Algorithms for choosing $\alpha_k$ often seek to satisfy conditions of sufficient decrease and curvature -- keeping $\alpha_k$ from growing too large or shrinking too small. Popular conditions are the Armijo conditions, Wolfe conditions, and Goldstein conditions with various relaxations of each \cite{Nocedal2006}[section 3], \cite{MoreThuente}, \cite{Bonnans2006}.

\textbf{Trust-region methods}

Trust-region methods \cite{Nocedal2006}[section 4] are another popular class of algorithms which seek to approximately solve the constrained optimisation problem
\begin{align}
    p^* &= \argmin_p m_k(p) \text{ subject to } \norm{p} \leq \Delta_k  \label{eqn:trust-region} \\
    x_{k + 1} &= x_k + p^*
\end{align}
for some norm $\norm{\cdot}$ and trust-region radius $\Delta_k \in \RR^+$. The trust-region radius is chosen at each step based upon how well $m_k$ approximated $f$ at the previous step, see \cite{Conn2000}[sections 6.1 and 10.5] for details.

The minimum in (\ref{eqn:trust-region}) is found using a variety of approximate methods \cite{Nocedal2006}[section 4],  \cite{steihaug}, and when $\norm{-H_k^{-1} \grad f(x_k)} \leq \norm{\Delta_k}$, then by construction $p^* = -H_k^{-1} \grad f(x_k)$. For this reason, trust-region methods are often interpreted as a class of methods for interpolating between $x_k$ and $-H^{-1} \grad f(x_k)$ via the scalar parameter $\Delta_k \in \left[0, \norm{-H^{-1} \grad f(x_k)}\right].$

One advantage line search algorithms have over trust-region algorithms is the value $H_k^{-1} \grad f(x_k)$ can be computed once, cached, and reused as $\alpha_k$ varies. This is not always true for trust-region algorithms, which may require significant effort to recompute $p^*$ as $\Delta_k$ varies.


\textbf{Searches in Optimistix}

A \emph{search} is a new abstraction introduced in Optimistix to generalise line searches, trust-region methods, and learning-rates. Searches are defined as functions taking local information and an internal state, and producing a scalar
\begin{equation}
    s: \mathcal{D} \times \mathcal{S} \to \RR.
\end{equation}
\noindent

For an example of a line search, consider the backtracking Armijo update. The Armijo backtracking algorithm uses local data $d_k = (d_k^{(1)}, d_k^{(2)}) \in \mathcal{D} = \RR \times \RR^N$, where $d_k^{(1)}$ represents the objective function value $f(x_k)$, and $d_k^{(2)}$ it's gradient $\grad f(x_k)$. The Armijo search state is $\sigma_k = (\sigma_k^{(1)}, \sigma_k^{(2)}) \in \mathcal{S} = (0, 1] \times \{0, 1\}$, where $\sigma_k^{(1)}$ represents the current step-size and $\sigma_k^{(2)}$ represents whether to shrink or reset the step size, depending on whether the last step satisfied the Armijo condition. For a decrease factor $c \in (0, 1]$, the Armijo search takes the step
\begin{equation}
    s(d_k, \sigma_k) = \begin{cases}
        c \sigma_k^{(1)} & \sigma_k^{(2)} = 0 \\
        1 & \sigma_k^{(2)} = 1, \label{eqn:armijo-search}
    \end{cases}
\end{equation}
and updates its state $\sigma_k$ via
\begin{equation*}
    \sigma_{k+1}^{(1)} = s(d_k, \sigma_k)
\end{equation*}
and $\sigma_{k+1}^{(2)} = 1$ when the Armijo condition
\begin{equation}
    f(x_k + \delta_k) \leq d_k^{(1)} + \eta\delta_k^T d_k^{(2)} \label{eqn:armijo}
\end{equation}
is satisfied, or $\sigma_{k + 1}^{(2)} = 0$ otherwise. Here, $0 < \eta < 1$ is a hyperparameter which determines how much a step must decrease to be accepted (larger means more decrease is required) and $\delta_k$ is the proposed step.

There are a number of different line search algorithms, but trust-region methods typically use the same
trust-region radius selection algorithm. This algorithm is represented by the search
\begin{equation}
    s_\text{TR}(d_k, \sigma_k) = \begin{cases}
        c_2 \sigma_k^{(1)} & \sigma_k^{(2)} > C_2, \\
        \sigma_k^{(1)} & C_1 < \sigma_k^{(2)} < C_2 \\
        c_1 \sigma_k^{(1)} &  \sigma_k^{(2)} < C_1 \\
    \end{cases}
\end{equation}
where $0 < c_1 < 1$ is a decrease amount, $c_2 > 1$ an increase amount,
$C_1$ a low cutoff (close to $0$) and $C_2$ is a high cutoff (close to $1$.)
The state is updated via
\begin{equation*}
    \sigma_{k+1}^{(1)} = s(d_k, \sigma_k)
\end{equation*}
and
\begin{equation}
    \sigma_{k}^{(2)} = \frac{f(x_k) - f(x_k + \delta_k)}{m_k(0) + m_k(\delta_k)}
\end{equation}
where $m_k$ is the quadratic model function (\ref{eqn:preconditioned}) and $\delta_k$ is again the
proposed step $x_{k + 1} = x_k + \delta_k$. The second component of the state, $\sigma_k^{(2)}$, is referred to as the `trust-region ratio', and roughly indicates how well
the model function $m_k$ predicted the decrease in $f$ that would come from taking the step $\delta_k$.

Current searches implemented in Optimistix are: learning rate, backtracking Armijo line search \cite{Nocedal2006}[section 3.1], the classical trust-region ratio update \cite{Conn2000}[section 6.1], and a trust-region update using a linear local approximation for first-order methods \cite{Conn2000}.

\subsection{Descent}

A \emph{descent} is another new abstraction in Optimistix. It is a function taking a scalar, the same local data, and an internal state and producing a vector:
\begin{equation}
    \delta: \RR \times \mathcal{D} \times \mathcal{R} \to \RR^N.
\end{equation}
\noindent
The value produced by $\delta$ is used to update the iterates $x_k$
\begin{equation*}
    x_{k + 1} = x_k + \delta(\alpha_k, d_k, \gamma_k)
\end{equation*}
\noindent
with $\alpha_k \in \RR, \ d_k \in \mathcal{D}, \ \gamma_k \in \mathcal{R}$.


The main function of the descent is to map the scalar produced by
the search $\alpha_k$ into a meaningful optimisation step.

For example, we can represent the preconditioned gradient step used in a
line search as a descent via
\begin{equation}
    \delta_\text{LS}(\alpha_k, d_k, \gamma_k) = -\alpha_k(d^{(2)})^{-1} (d^{(1)}) \label{eqn:pre_descent}
\end{equation}
\noindent
here, $d_k = (d^{(1)}_k, d^{(2)}_k) \in \mathcal{D} = \RR^N \times \RR^{N \times N}$.
Breaking this descent down, consider applying (\ref{eqn:pre_descent}) to a minimisation problem.
If $d_k = (\grad f(x_k), H_k)$, then (\ref{eqn:pre_descent}) is equivalent to 
the standard preconditioned gradient descent $\delta_\text{TR}(\alpha_k, d_k, \gamma_k) = -\alpha_k H_k^{-1} \grad f(x_k)$, 
with a step-size $\alpha_k$. If $d_k = (r_k, J_k)$, the residual vector $r_k$ and Jacobian $J_k$ at step $k$,
then (\ref{eqn:pre_descent}) is the Gauss-Newton step with step-size $\alpha_k$ for a least-squares problem.

As another example, we can represent the trust-region subproblem as a descent via
\begin{equation}
    \delta_\text{TR}(\alpha_k, d_k, \gamma_k) \approx
    \argmin_{p, \ \Vert p \Vert \leq \alpha_k} d_k^{(1)}+ \grad (d_k^{(2)})^T p + \frac{1}{2} p^T d_k^{(3)} p \label{eqn:tr-descent}
\end{equation}
for $d_k \in \mathcal{D} = \RR \times \RR^N \times \RR^{N \times N}$. For a minimisation
problem this is $d_k = (f(x_k), \grad f(x_k), H_k)$. Here, the output of the search
is mapped not to a line search, but to the trust-region radius of a trust-region
algorithm.

$\alpha_k$ has a similar meaning in both of these algorithms. It's used to restrict
the size of the update at step $k$ when confidence in the approximation $m_k$ is low.
By abstracting the map from $\alpha_k$ to the optimiser update $\delta_k$, descents allow
a search to be used in a variety of different algorithms without knowing exactly how 
the step-size will affect the optimiser update. This decoupling is novel to optimistix,
and is one of its most powerful features.

Optimistix descents, like classical line searches, allow certain values used in their computation to be cached and reused. For example, in $\delta_\text{LS}$, the value $(d^{(2)})^{-1} d^{(1)}$ is only computed once as $\alpha_k$ varies.






Current descents implemented in Optimistix are: steepest descent, nonlinear conjugate gradient descent \cite{Shewchuk94}, direct and indirect damped Newton (used in the Levenberg-Marquardt algorithm \cite{MoreLM, Nocedal2006}), and dogleg \cite{Nocedal2006}[section 4.1]. The indirect damped Newton and dogleg descents both approximately solve (\ref{eqn:trust-region}) (see \cite{Nocedal2006}[sections 4.1 \& 4.3]) and are capable of using any user-provided function which maps a PyTree to a scalar as the norm $\Vert \cdot \Vert$.




\textbf{Flattening the bilevel optimisation problem}

The classical approach to implementing a line search is as a bilevel optimisation problem. The inner loop to obtain the solution to the line search, the outer loop over the iterates $x_k$.

The termination condition of the line-search in the inner loop, such as the Armijo or Wolfe conditions, may evaluate $f(x_{k + 1})$ when determining whether to halt the inner loop. This can result in an extra function evaluation or additional logic at each step of the outer loop, as $f(x_{k + 1})$ is usually computed in each step of the outer loop as well (eg. when computing the gradient of $f$.)

Not every search has termination conditions which require computing $f(x_{k + 1})$, so we cannot avoid an extra function evaluation by requiring the inner loop to evaluate $f(x_{k + 1})$ and pass this to the outer loop.

Optimistix gets around this by using step-rejection: a step $x_{k + 1}$ is proposed, and the outer loop computes $f(x_{k+1})$, and at the next step the search decides whether to accept or reject $x_{k + 1}$.

\section{Converting Between Problem Types} \label{sec:interoperability}

Optimistix handles four types of nonlinear optimisation task, each with a simple user API:
\begin{itemize}
    \itemsep0em
    \item \begin{flushleft}
        Minimisation:
    \end{flushleft}
    \verb|optimistix.minimise|.
    \item \begin{flushleft}
        Nonlinear least-squares solves: 
    \end{flushleft} \verb|optimistix.least_squares|.
    \item \begin{flushleft}
        Root-finding:
    \end{flushleft} \verb|optimistix.root_find|.
    \item \begin{flushleft}
        Fixed-point iteration:
    \end{flushleft} \verb|optimistix.fixed_point|.
\end{itemize}
\noindent
For example, to find a root of the function \verb|fn| we call:
\begin{lstlisting}[language=python, numbers=none]
solver = optimistix.Newton(
  rtol=1e-3, atol=1e-3
)
optimistix.root_find(fn, solver, y0)
\end{lstlisting}
\noindent
where \verb|y0| is an initial guess and \verb|optimistix.Newton| is the optimistix solver used to find the root. \verb|rtol| and \verb|atol| are relative and absolute tolerances used to deteremine whether the solver should terminate, as described in section \ref{sec:termination}.

Optimistix can convert between these problem types when appropriate. A fixed-point iteration $f(x) = x$ is converted to a root-find problem by solving $f(x) - x = 0$. A root-find problem $f(x) = 0$ is converted to a nonlinear least-squares problem by using entries of $f(x)$ as the residuals: $\min_x \sum_i f(x)_i^2$. A least-squares problem is already a special case of minimisation problem, so we take the objective function $g(x) = \sum_i f(x)_i^2$ and provide this directly to a minimiser.

Converting between problem types is not just a neat trick, but a common pattern in optimisation, especially for root-finding and fixed-point iteration. While standard methods for root-finding and fixed-point iteration, such as Newton, chord, bisection, and fixed-point iteration \cite{Nocedal2006, Bonnans2006} are available in Optimistix, these methods tend not to work as well for complex, highly nonlinear problems. In this case, common practice is to convert the fixed-point or root-find problem into a least-squares or minimisation problem, and solve using a minimiser or nonlinear least-squares solver \cite{Nocedal2006}[section 11].

The conversion between problem types outlined above is done automatically: the solver is checked, and the problem automatically converted and lowered to a solve of that type. For example,
\begin{lstlisting}[language=Python, numbers=none]
solver = optimistix.BFGS(
  rtol=1e-3, atol=1e-3
)
optimistix.root_find(fn, solver, y0)
\end{lstlisting}
will convert the root-find problem to a minimisation problem and use the BFGS algorithm \cite{Nocedal2006}[section 6.1] to find a root.


\section{Automatic Differentiation}
\label{sec:diff}

All solvers in Optimistix are iterative, and provide two methods for automatic differentiation:
\begin{itemize}
    \item Differentiation via the implicit function theorem \cite{jaxopt_implicit_diff}
    \item Differentiation via online treeverse \cite{stumm2010new, wang2009minimal}.
\end{itemize}
The former is the default for both forward-mode autodiff and backpropagation in Optimistix. Both are known automatic differentiation techniques; however, to the best of our knowledge, Optimistix is the first nonlinear optimisation library in JAX to feature online treeverse. We now discuss each of these in turn.

\subsection{Forward-Autodiff and Backpropagation Using the Implicit Function Theorem}

Consider a continuously differentiable function $F:\RR^n \times \RR^m \to \RR$ and values $x_0 \in \RR^n$ and $\theta_0 \in \RR^m$ satisfying
\begin{align}
    F(x_0, \theta_0) &= 0 \label{eqn:IFTassumptions1}\\
    \det\left(\frac{\D F}{\D x}(x_0, \theta_0)\right) &\neq 0. \label{eqn:IFTassumptions2}
\end{align}
The implicit function theorem (IFT) states there exists neighborhoods $x_0 \in N_{x_0}$ and $\theta_0 \in N_{\theta_0}$ and a differentiable function $x^*: N_{\theta_0} \to N_{x_0}$ such that
\begin{align*}
    x^*(\theta_0) &= x_0 \\
    F(x^*(\theta), \theta) &= 0 \ \ \forall x \in N_{\theta_0}.
\end{align*}
and 
\begin{equation}
    \frac{\D x^*}{\D \theta} = - \left(\frac{\D F}{\D x} (x_0, \theta_0)\right)^{-1} \frac{\D F}{\D \theta} (x_0, \theta_0) \label{eqn:IFTdiff}
\end{equation}
The implicit function theorem gives means of calculating the derivative $\frac{\D x^*}{\D\theta}$, despite the fact that $x^*$ is only defined implicitly, and is the solution to a nonlinear problem.

Every optimisation task in Optimistix (minimisation, ...) may be rewritten in the form of equation (\ref{eqn:IFTassumptions1}), and thus the derivative of the solution (the argminimum, ...) with respect to any parameters $\theta$ may be found using equation (\ref{eqn:IFTdiff}).

For root finding, equations (\ref{eqn:IFTassumptions1}-\ref{eqn:IFTassumptions2}) apply directly. For fixed-point iteration, the system of nonlinear equations is $F(x, \theta) = f(x, \theta) - x$. For minimisation and least-squares, $F(x, \theta) = \frac{\D f}{\D x}(x,\theta)$, as the minimum is found at a critical points of $f$.

We highly recommend \cite{jaxopt_implicit_diff} for many more details on using the implicit function theorem to differentiate nonlinear optimisation routines.

\subsection{Backpropagation using Online Treeverse}

An optimisation routine can also be differentiated naively, by backpropagating through the components of the optimiser at each iteration $k$. This method works even when the assumptions (\ref{eqn:IFTassumptions1}-\ref{eqn:IFTassumptions2}) of the implicit function theorem are not satisfied, or the full Jacobian $\frac{dF}{dx}$ cannot be constructed due to memory constraints. Further, it can be faster in some cases, especially when the dimensionality of the domain of the optimisation problem is much larger than the number of iterations $K$. See \cite{ablin2020superefficiency} for more details comparing the efficiency of these methods.

The main downside of the naive backpropagation algorithm for nonlinear optimisation is memory cost. It costs $\mathcal{O}(KM)$ memory to backpropagate through an optimisation, were $M$ is the memory cost of backpropagating through a single iteration of the optimiser, and $K$ is the total number of iterations used. The classical way around this limitation is checkpointing, which chooses a fixed number $N \in \mathbb{N}$ of checkpoint iterations, and stores the values of the optimiser only at those checkpointed iterations, reducing the memory cost to $\mathcal{O}(NM)$. Iterations which are not checkpointed are recomputed from the previous checkpoint, trading off memory cost for extra computation time.

Online treeverse is a checkpointing backpropagation scheme designed for performing backpropagation over loops where the number of iterations is a-priori unknown. Online treeverse dynamically updates its checkpoint locations as the loop continues, keeping optimal computation times for a fixed number of checkpoints $c$. For further details see  \cite{stumm2010new} and \cite{wang2009minimal}.

\section{Termination Condition} \label{sec:termination}



Optimistix introduces a Cauchy convergence criteria, where iterations are stopped if 
\begin{equation*}
    \abs{\frac{f(x_{k + 1}) - f(x_{k})}{\eps_a + \eps_r f(x_k)}} < 1 \text{ and } \norm{\frac{x_{k + 1} - x_{k}}{\eps_a + \eps_r x_k}} < 1
\end{equation*}
for $\eps_a$ an absolute tolerance set by the user and $\eps_r$ a relative tolerance set by the user, and division is performed elementwise. By default, the norm $\norm{\cdot}$ is $\norm{\cdot}_\infty$, ie. the maximum over the elementwise absolute differences in $x_{k + 1}$ and $x_k$. However, users can set this themselves to be any function from a PyTree to a scalar.

To justify this choice, we note that the optimisation literature is not consistent with the termination criteria used \cite{Nocedal2006}, \cite{Bonnans2006}, \cite{Conn2000}, \cite{minpack}. For example, \cite{sasuserguide} outlines 20 different convergence criteria currently in use for nonlinear optimisation.

Without an established convention, we choose this termination condition to match the better established convention for the literature on solving numerical differential equations \cite{hairer2008solving-i, hairer2002solving-ii}. This is reasonable choice for optimisation, resembling (but not exactly matching) a combination of the `X-convergence' and `F-convergence' criteria in the highly successful MINPACK \cite{minpack} optimisation suite. Alternatively, it is roughly a combination of all four of the `FTOL', `ABSFTOL', `XTOL', and `ABSXTOL' criteria found in \cite{sasuserguide}.  This also ensures a consistent approach between Optimistix and existing JAX + Equinox libraries, such as Diffrax \cite{KidgerThesis}.

\section{Experiments}
\label{sec:experiments}

\subsection{Summary of Experiments}
We demonstrate the fast compilation and run times of Optimistix by benchmarking our 
solvers for both runtime and compile time against the state-of-the-art optimisation algorithms in JAXopt \cite{jaxopt_implicit_diff} and SciPy \cite{scipy}. For SciPy, we only benchmark runtimes, since SciPy does not have a notion of compile times (however, we still compile the function passed to the SciPy solver. The overhead of this function compilation is excluded.)

\begin{table}[]
    \centering
    \caption{Average minimum runtime comparison}
    \begin{tabular}{ |p{2.1cm}||p{1cm}|p{1cm}|p{1cm}|  }
        \hline
        \multicolumn{4}{|c|}{Average minimum runtime (milliseconds)} \\
        \hline
        Optimiser & Optx & JAXopt & SciPy\\
        \hline
        BFGS & \textbf{0.836} & 16.4 & 250 \\
        Nonlinear CG & \textbf{0.466}  & 6.27 & 208\\
        LM1 & - & \textbf{1.70} & -\\
        LM2 & \textbf{144} & - & 188\\
        Gauss-Newton & 4.54 & \textbf{3.08} & N/A\\
        \hline
    \end{tabular} 
    \label{tab:runtime}
\end{table}

\begin{table}[]
    \centering
    \caption{Average compile time comparison}
    \begin{tabular}{ |p{2.1cm}||p{2cm}|p{2cm}|  }
        \hline
        \multicolumn{3}{|c|}{Average compile time (milliseconds)} \\
        \hline
        Optimiser & Optx &JAXopt\\
        \hline
        BFGS & \textbf{251} & 854 \\
        Nonlinear CG & \textbf{130}  & 849\\
        LM1 & - & \textbf{554}\\
        LM2 & \textbf{296} & -\\
        Gauss-Newton &\textbf{288} & 416\\
        \hline
    \end{tabular}
    \label{tab:compile_time}
\end{table}

Tables \ref{tab:runtime}-\ref{tab:compile_time} present the average minimum runtime and average compile time for a set of over 100 test problems. Each problem is solved 10 times to a fixed tolerance with a max budget of 2000 iterations, and the minimum over those 10 iterations is taken to be the runtime of the solver for that problem. The reported values are averaged over all the problems in the test set. The full methodology and additional comparisons are described in appendix \ref{appendix:comparison}.

The average compile times for Optimistix (Optx) are up to 6.5 times faster than the average compile times for JAXopt. The runtimes include all operations for a full solve; notably, it includes the difference in the number of iterations arising from convergence criteria and step-sizes in Optimistix and JAXopt. The times for a single iteration are presented in appendix \ref{appendix:comparison}, where JAXopt is marginally faster for a single iteration than Optimistix.

The difference between Levenberg-Marquardt (LM) algorithms LM1 and LM2 relates to a technical detail in their implementation. At each iteration, LM performs a linear solve of the equation $(J_k^TJ_k + \lambda I)p = -J^T r_k$, where $J_k$ is the Jacobian and $r_k$ is the residual of the objective function at step $k$. This linear system can be solved directly with a Cholesky or conjugate gradient solver, but the condition number of $J_k$ is squared in $J_k^TJ_k$. This is the approach used in LM1, which JAXopt implements. Alternatively, this linear system can be transformed into an equivalent least-squares problem which does not square the condition number, but which requires a slower linear solve. This is the approach used in LM2, which Optimistix and SciPy use. We accept the slower solve in return for increased robustness and accuracy. This is detailed in \cite{MoreLM} and \cite{Nocedal2006}[section 10]. A similar tradeoff is made for Gauss-Newton. The difference in robustness is noticeable: Optimistix LM failed to solve 3 of the test problems to the specified accuracy within 2000 iterations, SciPy failed to solve 7, and JAXopt failed to solve 13.

Compile times in the table are 3-4 orders of magnitude longer than runtimes. For solves with fewer than many thousands of iterations, the total time is dominated by compilation, where Optimistix is faster.

\subsection{Benchmarking with Performance Profiles}
Performance profiles are an established technique for benchmarking
optimisation software \cite{performanceprofiles, Beiranvand_2017}. Performance profiles are defined in terms of solver performance ratios.




Formally, let $\mathcal{S}$ be a collection of optimisers, and $\mathcal{P}$ a collection of problems. Letting
\begin{align*}
    R_{s, p} &= \text{the time for optimiser } s \text{ to solve problem } p, \\
    &\hspace{0.5cm} \text{excluding compilation time.}\\
    C_{s, p} &= \text{the time for optimiser } s \text{ to compile problem } p,
\end{align*}
the runtime performance ratio and compile time performance ratios are:
\begin{align*}
    r^R_{s, p} &= \frac{R_{s, p}}{\min\{R_{s, p}\colon s \in \mathcal{S}\}} \\
    r^C_{s, p} &= \frac{C_{s, p}}{\min\{C_{s, p}\colon s \in \mathcal{S}\}}
\end{align*}
The runtime and compile time performance profiles are then:
\begin{align*}
    \rho^R_s(\tau) = \frac{1}{\vert\mathcal{P}\vert} |\{p \in \mathcal{P}\colon r^R_{s, p} \leq \tau \}| \\
    \rho^C_s(\tau) = \frac{1}{\vert\mathcal{P}\vert} |\{p \in \mathcal{P}\colon r^C_{s, p} \leq \tau \}|
\end{align*}
ie. $\rho^R_s(\tau)$ is the proportion of problems solver $s$ solved within a
factor $\tau$ of the minimum runtime attained by any solver in $\mathcal{S}$.
The value $\rho^R_s(1)$ is the proportion of problems where solver $s$ had the best
runtime performance  in $\mathcal{S}$. For large values of $\tau$, $\rho^R_s(\tau)$ represents the proportion of problems solver $s$ accurately solved. Similarly, $\rho^C_s(\tau)$ is the proportion of problems
solved within a factor of $\tau$ of the best compile time.

\begin{figure}
    \centering
    \includegraphics[width=\linewidth]{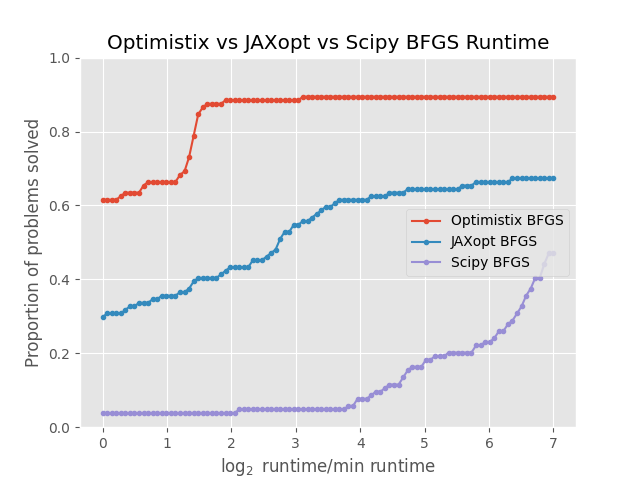}
    \caption{BFGS runtime comparison}
    \label{fig:allbfgs}
\end{figure}

\begin{figure}
    \centering
    \includegraphics[width=\linewidth]{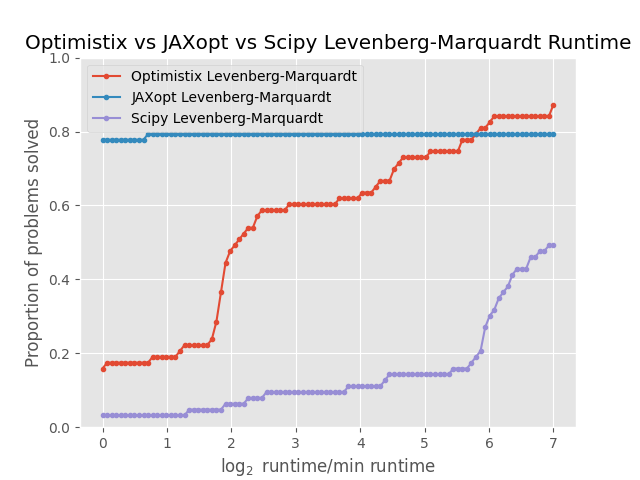}
    \caption{Levenberg-Marquardt runtime comparison}
    \label{fig:allLM}
\end{figure}

\begin{figure}
    \centering
    \includegraphics[width=\linewidth]{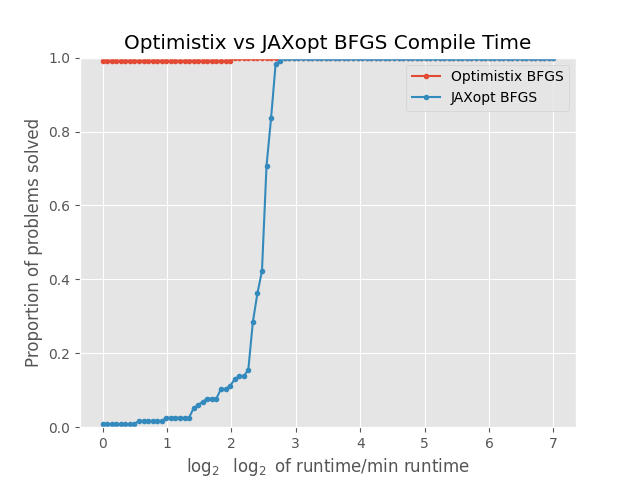}
    \caption{BFGS compile time comparison}
    \label{fig:bfgscompile}
\end{figure}

\begin{figure}
    \centering
    \includegraphics[scale=0.5]{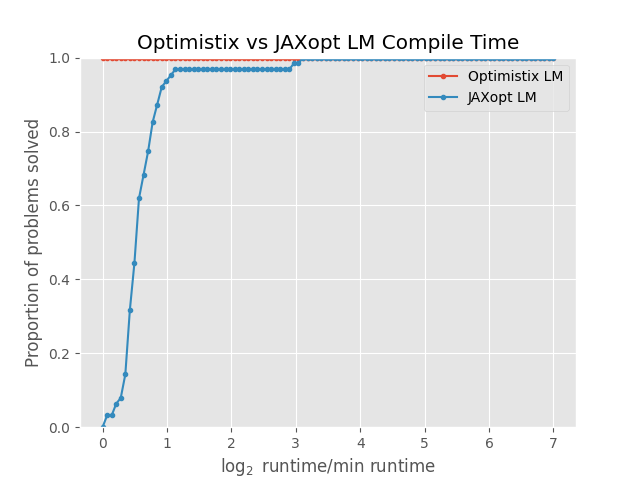}
    \caption{Levenberg-Marquardt compile time comparison}
    \label{fig:LMcompile}
\end{figure}

In figures \ref{fig:allbfgs}-\ref{fig:allLM} we see the runtime performance profiles of BFGS, and Levenberg-Marquardt with Optimsitix in red, JAXopt in blue, and SciPy in purple (no distinction is made between LM1 and LM2 in the plot.) The plots are log scale on the x-axis (representing $\tau$ in the performance profile,) and range from $1$ to $2^7$. Optimistix performs the best of all solvers for BFGS. For Levenberg-Marquardt it is generally slower than JAXopt but solves more problems in total.

It may appear as though SciPy has not solved many of the test problems; however, it is just more than 250 times slower and therefore does not show on the plots. 

In figures \ref{fig:bfgscompile}-\ref{fig:LMcompile} we see a compile time comparison of BFGS and Levenberg-Marquardt. Optimistix outperforms JAXopt on every problem for both optimisers.

\section{Related Work}

\textbf{Gradient-based minimisation    } Optax \cite{deepmind2020jax} implements algorithms for gradient-based minimisation in JAX. Optimistix has a broader scope than Optax, handling other nonlinear optimisation tasks and second-order optimisers. Optax and Optimistix are compatible libraries, and Optax minimisers can be used within Optimistix via \verb|optimistix.OptaxMinimiser|.

\vspace{0.5cm}
\noindent
\textbf{General-purpose optimisation    } JAXopt \cite{jaxopt_implicit_diff} is an excellent differentiable linear and nonlinear optimisation library in JAX which includes many optimisation tasks, including ones out-of-scope for Optimistix such as quadratic programming and non-smooth optimisation.

We see Optimistix and JAXopt as having fundamentally different scopes and core abstractions. JAXopt has a larger scope than Optimistix, but its focus is not on modularity. While specific choices of line search algorithms can be interchanged in some cases, for the most part introducing a new optimiser requires writing the algorithm in its entirety. We recommend JAXopt for general optimisation tasks.

SciPy \cite{scipy} also offers general-purpose nonlinear optimisation routines, many of which call into MINPACK \cite{minpack}, a nonlinear optimisation library written in Fortran. SciPy includes many optimisation tasks, including minimisation, least-squares, root-finding, and global optimisation. However, SciPy implementations are not differentiable, and are generally difficult to extend.

\section{Conclusion}

We introduced Optimistix, a nonlinear optimisation software for minimisation, least-squares, root-finding, and fixed-point iteration with a novel, modular approach.

\section{Impact Statement}
This paper presents work whose goal is to advance the field of Machine Learning. There are many potential societal consequences of our work, none which we feel must be specifically highlighted here.

\bibliographystyle{icml2023}
\bibliography{main}


\newpage
\appendix
\onecolumn
\section{Experiment details}
\label{appendix:comparison}

\subsection{Methodology}

We choose as our comparison set BFGS, nonlinear CG, Levenberg-Marquardt, and
Gauss-Newton, as these are the four nontrivial minimisation/least-squares algorithms shared
by Optimistix and JAXopt (ie. gradient descent is excluded intentionally.) Where solver implementations differ, we try and match them as closely as possible.

Specifically, where algorithms differ from the base implementations of the algorithms is:
\begin{itemize}
    \item BFGS uses backtracking Armijo line search in both implementations (default for JAXopt is zoom.)
    \item Gauss-Newton uses a conjugate gradient solver on the normal equations in both implementations (default for Optimisix is the polyalgorithm \verb|AutoLinearSolver(well_posed=None)| from \cite{lineax}.
\end{itemize}

The set of test problems consists of 104 minimisation problems, 63 of which are least-squares problems, taken from the test collections \cite{Jamil_2013}, \cite{numericaleval}, \cite{Andrei2008AnUO}, \cite{minpack2}, and \cite{tuos}. Levenberg-Marquardt and Gauss-Newton are only ran on the least-squares problems. Solvers are initialised at canonical initialisations when available (see \cite{tuos} and \cite{Andrei2008AnUO}).

Runtime is a noisy measurement. During an experiment, the computer an experiment is running on may be running background processes which we cannot control. To mitigate this, when assessing runtime we run each problem 10 times and take the minimum over these repeats. This indicates roughly what the "best we can expect" from a given solver is.

The global minima are not known for all test problems. On problems where minima are known, we require that 
\begin{equation*}
    \frac{|f(x_N) - f(x^*)|}{\epsilon_a + \epsilon_r f(x^*)} < 1
\end{equation*}
for an absolute tolerance $\epsilon_a$ and relative tolerance $\epsilon_r$, and $x_N$ the argmin found by the solver. If this condition fails, the runtime is set to \verb|jnp.inf|. This information is automatically incorporated in to the performance profile; however, it is not included in the tables of average runtime performance. This is because there is no obvious way to penalise nonconvergence when comparing runtimes. To get around this, we assure that in all cases the number of failures for an Optimistix solver was less than or equal to the number of failures for a JAXopt solver, as to not give an unfair advantage to Optimistix.

Though this convergence criteria looks similar to the Cauchy termination condition for Optimistix, it is applied to $f(x_N) - f(x^*)$ and not $f(x_{k+1}) - f(x_{k})$. $f(x^*)$ is an unknown quantity to both solvers, so there is no advantage provided by the convergence criterion in Optimistix.

Finally, while the general formulation of performance profile allows for $\mathcal{S}$ to contain more than two optimisers, using performance profiles with more than two optimisers requires careful interpretation. For example, when comparing three optimisers, it is possible that the method which appears to be second best performs worse than the method which appears to be third best when these two methods are compared directly.

For this reason, we provide all the pairwise comparisons of Optimistix vs JAXopt and SciPy in \ref{sec:further_experiments}, with the Optimsitix solver in red and the comparison solver (JAXopt or SciPy) in blue.

\subsection{Further experiments}
\label{sec:further_experiments}

For simplicity, we included a number of performance profiles with one-on-one comparisons of Optimistix vs JAXopt and Optimistix vs SciPy. Also included is a table of average minimum runtimes for a single optimiser run. Ultimately, no user will use only a single optimiser run, so the aggregate information presented in section \ref{sec:experiments} should be taken as more informative.


\newpage

\centering
\begin{tabular}{ |p{3cm}||p{3cm}|p{3cm}|p{3cm}|  }
    \hline
    \multicolumn{4}{|c|}{Average minimum runtime for a single optimiser run (milliseconds)} \\
    \hline
    Optimiser & Optimistix &JAXopt&SciPy\\
    \hline
    BFGS & 0.151 & \textbf{0.138} & 2.20 \\
    Nonlinear CG & 0.123  & \textbf{0.0645}   & 2.61\\
    LM & 0.515 & \textbf{0.0917} &  3.31\\
    Gauss-Newton &0.393 & \textbf{0.0578} &  N/A\\
    \hline
\end{tabular} \label{table:single_run}

\begin{figure}[h]
    \begin{subfigure}
        \centering
        \includegraphics[scale=0.5]{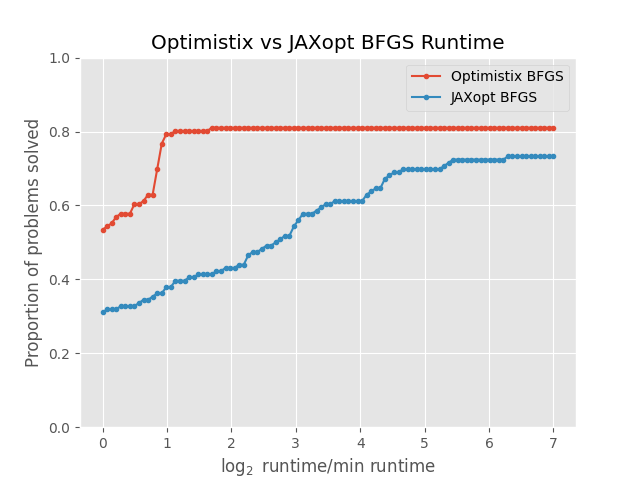}
        \label{fig:bfgs}
    \end{subfigure} \hfill
    \begin{subfigure}
        \centering
        \includegraphics[scale=0.5]{bfgs_compile_times.png}
    \end{subfigure} \hfill
    \caption{Optimistix vs JAXopt BFGS}
\end{figure}

\begin{figure}[h]
	\begin{subfigure}
		\centering
		\includegraphics[scale=0.5]{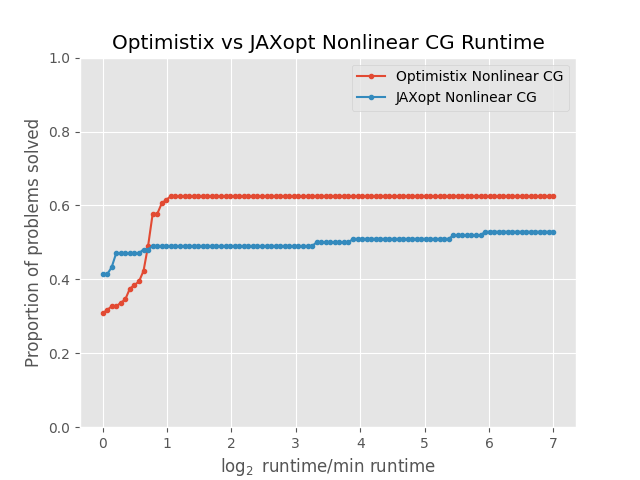}
	\end{subfigure}\hfill
	\begin{subfigure}
		\centering
		\includegraphics[scale=0.5]{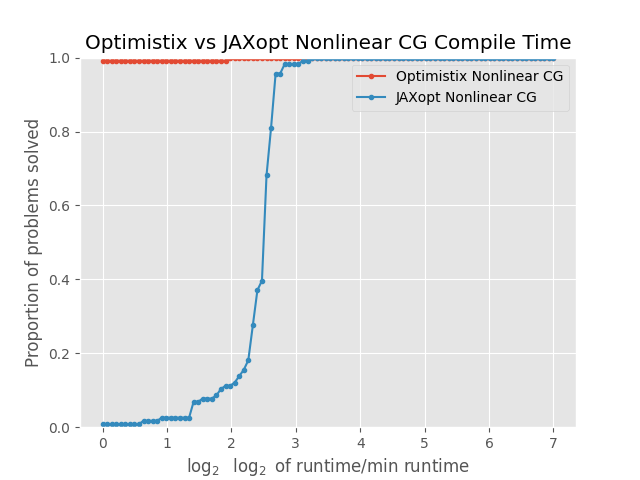}
	\end{subfigure} \hfill
        \caption{Optimistix vs JAXopt nonlinear CG}
\end{figure}

\begin{figure}[h]
	\begin{subfigure}
		\centering
		\includegraphics[scale=0.5]{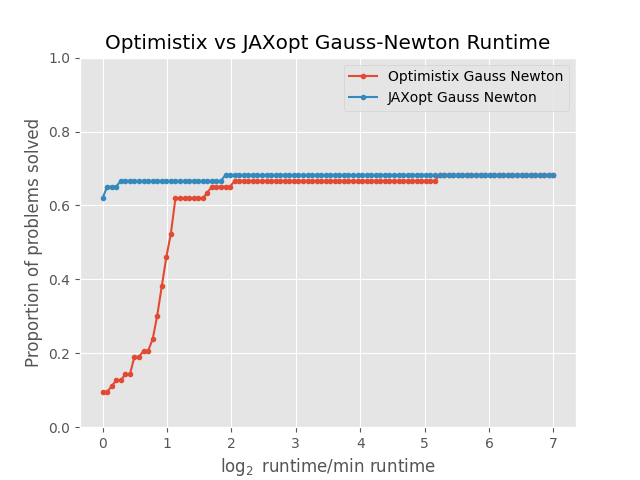}
	\end{subfigure}\hfill
	\begin{subfigure}
		\centering
		\includegraphics[scale=0.5]{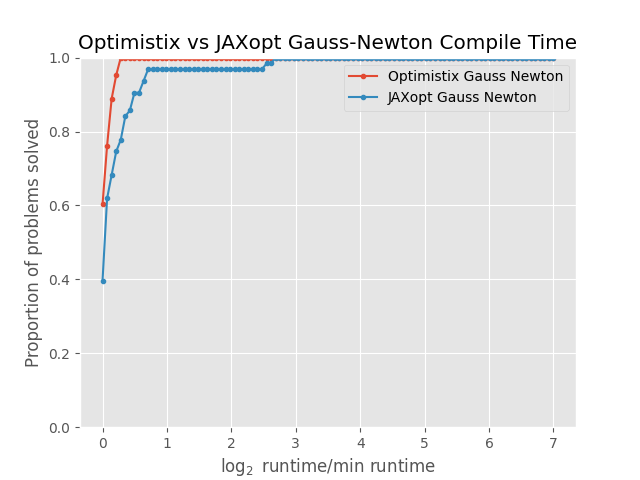}
	\end{subfigure}
	\caption{Optimistix vs JAXopt Gauss-Newton}
	\label{fig:gn}
\end{figure}

\begin{figure}[h]
	\begin{subfigure}
		\centering
		\includegraphics[scale=0.5]{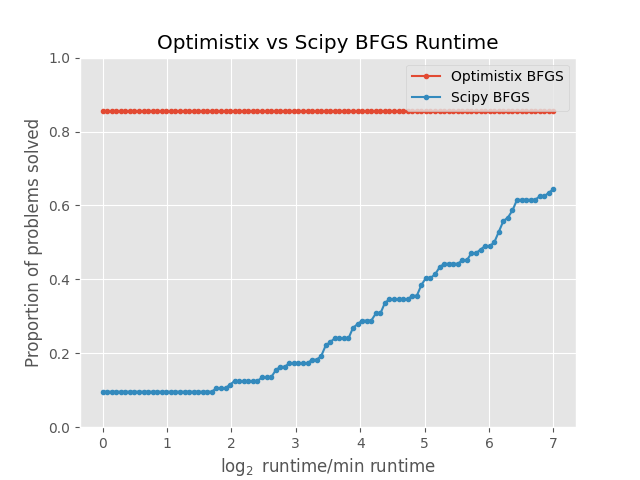}
	\end{subfigure}\hfill
	\begin{subfigure}
		\centering
		\includegraphics[scale=0.5]{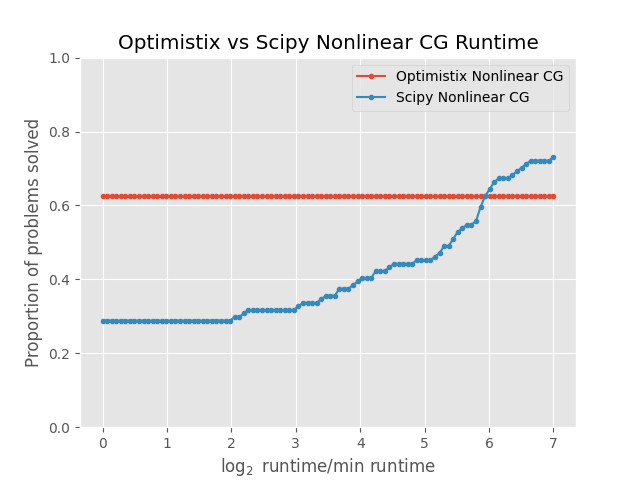}
	\end{subfigure}
	\caption{Optimistix vs SciPy BFGS (left) and nonlinear CG (right)}
	\label{fig:gn}
\end{figure}

\begin{figure}
    \centering
    \includegraphics[scale=0.5]{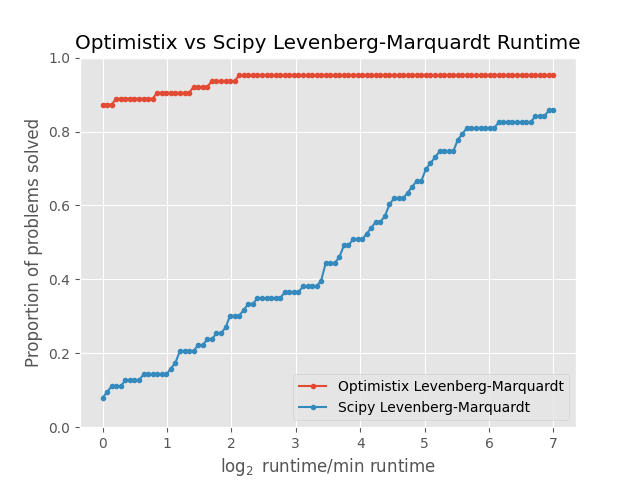}
    \caption{Optimistix vs SciPy Levenberg-Marquardt runtimes}
    \label{fig:scipylm}
\end{figure}


\end{document}